\documentclass{amsart}
\usepackage{amssymb, amsmath, latexsym, amsthm}
\usepackage{fancyhdr}
\usepackage{graphicx}

\newtheoremstyle{dotless}{}{}{\itshape}{}{\bfseries}{}{ }{}
\theoremstyle{dotless}

\theoremstyle{definition}
\DeclareMathOperator{\supp}{supp}

\DeclareMathOperator{\Lip}{Lip}
\def\supp{\operatorname{supp}}

\newcommand{\Om}{\Omega}

\newcommand{\cP}{\mathcal P}

\newcommand{\cS}{\mathcal S}

\begin{document}

\title[Riesz bases and frames 
   measured in the Kantorovich--Rubinstein norm]%
{Sign intermixing for Riesz bases and frames measured in the Kantorovich--Rubinstein norm}
\author[N.~Nikolski]{Nikolai Nikolski}
\address{Department of Mathematics,Universit\`{e} de Bordeaux, France}
\email{Nikolai.Nikolski@math.u-bordeaux.fr \textrm{(N.\ Nikolski)}}
\thanks{NN is partially supported by Grant MON 075-15-2019-1620 of the Euler International Mathematical  
Institute, St.Petersburg}
\author[A.~Volberg]{Alexander Volberg}
\thanks{AV is partially supported by the NSF  DMS 1900268 and by Alexander von Humboldt foundation}
\address{Department of Mathematics, Michigan Sate University, East Lansing, MI. 48823 USA and Hausdorff Center for Mathematics, Bonn 53115, Germany}
\email{volberg@math.msu.edu \textrm{(A.\ Volberg)}}
\makeatletter
\@namedef{subjclassname@10010}{
  \textup{10010} Mathematics Subject Classification}
\makeatother
% 42B	Harmonic analysis in several variables
% 42B100	Singular and oscillatory integrals (Calder?on-Zygmund, etc.)
% 42B35	Function spaces arising in harmonic analysis
% 47A	General theory of linear operators
% 47A100	Norms (inequalities, more than one norm, etc.)
%{100E100, 47B37, 47B40, 100D55.}
%
% 100D55	$H^p$-classes (1980-10009)
% 100E100	Integration, integrals of Cauchy type, integral representations of analytic functions
%
% 47B   	Special classes of linear operators
% 47B37	Operators on special spaces (weighted shifts, operators on sequence spaces, etc.)
% 47B40	Spectral operators, decomposable operators, well-bounded operators, etc.
\begin{abstract} 
We measure a sign interlacing phenomenon for Bessel 
sequences $ (u_{k})$ in $ L^{2}$ spaces in terms of the Kantorovich--Rubinstein 
mass moving norm $ \Vert u_{k}\Vert _{KR}$. Our main observation shows 
that, quantitatively, the rate of the decreasing $ \Vert u_{k}\Vert _{KR}\longrightarrow 
0$ heavily depends on S. Bernstein $ n$-widths of a compact of Lipschitz 
functions. In particular, it depends on the dimension of the measure 
space. 
\end{abstract}
\maketitle

\ \par 
  \section{What this note is about.}
\ \par 
 Let $ (\Omega ,\rho )$ be a metric space, and $ 
m$ a finite continuous (with no point masses) Borel measure on $ \Omega $. 
It is known \cite{NV2019} that for every \it frame $ (u_{k})_{k\geq 1}$ 
in $ L^{2}_{{\Bbb R}}(\Omega ,m)$\rm , the ``$ l^{2}$-masses" of positive 
and negative values $ u_{k}^{\pm }(x)$ are infinite:\ \par 
\ \par 
  \centerline{ $ \sum _{k}u_{k}^{+}(x)^{2}=$ $ \sum _{k}u_{k}^{-}(x)^{2}=$ 
$ \infty $ a.e. on $ \Omega $}
\ \par 
\rm (and moreover, $ \forall f\in L^{2}_{{\Bbb R}}(\Omega )$, $ f\not= 0$ 
$ \Rightarrow $ $ \sum _{k}(f,u_{k}^{\pm })^{2}_{L^{2}}=\infty $), 
where as usual $ u_{k}^{\pm }(x)=$ $ max(0,\pm u_{k}(x))$, $ x\in (0,1)$. 
So, at almost every point $ x\in \Omega $, there are many positive 
and many negative values $ u_{k}(x)$. Here, we show that for a fixed 
$ k$, positive and negative values are heavily intermixed.\ \par 
Precisely, we show that the measures $ u_{k}^{\pm 
}dm$ should be closely interlaced, in the sense that the Kantorovich-Rubinstein 
(KR) distances $\Vert u_{k}\Vert _{KR}=\Vert u_{k}^{+}-u_{k}^{-}\Vert 
_{KR}$ (see below) must be small enough. It is easy to see that if 
the supports $ \supp(u_{k}^{\pm })$ are distance separated from each 
other than $\Vert u_{k}\Vert _{KR}\approx\Vert u_{k}\Vert 
_{L^{1}  (m)}$, whereas in reality, as we will see, these norms are 
much smaller, and so, the sets $ \{x:u^{+}(x)>0\}$ and $ \{x:u^{+}(x)<0\}$ 
should be increasingly mixed. In this connection, it is interesting 
to recall one of the first (and classical) results in this direction, 
that of O. Kellogg \cite{Ke1916}, showing that on the unit interval $ \Omega =
I=:(0,1)$, the consecutive supports $ \supp(u_{k}^{\pm })$ are 
interlacing under quite general hypothesis on an orthonormal sequence 
$ (u_{k})$. (Later on, the sign interlacing phenomena were intensively 
studied for (orthogonal) polynomials (starting from P. Chebyshev, and 
earlier, see any book on orthogonal polynomials), so that, quite a 
recent survey \cite{Fi2008} counts about 780 pages and hundreds references; 
many new quantitative results are also presented).\ \par

\bigskip

Our results are most complete for the classical case 
$ \Omega =\, I^{d}$ ($ d\geq 1$) in $ {\Bbb R}^{d}$, $ I=(0,1)$, 
and $ m=\, m_{d}$ the Lebesgue measure and $ \rho $ the Euclidean 
distance on the cube. They also suggest that in general, the magnitudes 
of $ \Big \Vert u_{k}\Big \Vert _{KR}$ are defined by certain (unknown) 
interrelations between $ m$ and $ \rho $, and by a kind of the dimension 
of $ \Omega $. In fact, all depends on and is expressed in terms of 
a compact subset $ Lip_{1}$ of Lipschitz functions in $ L^{2}(\Omega ,m)$.\ \par 
\, \, \, Plan of the rest:\ \par 
2. Definitions and comments\ \par 
3. Statements on the generic behaviour of $ \Big \Vert u_{k}\Big \Vert _{KR}$\ \par 
4. Proofs\ \par 
5. Further examples and comments; numerical examples to Theorem 3.2; 
direct comparisons $ \Big \Vert u_{k}\Big \Vert _{KR}$ with Bernstein 
widths $ b_{k}(Lip_{1})$; an explicit expression for $ \Big \Vert u\Big \Vert 
_{KR}$.\ \par 
6. The fastest rates of decreasing $ \Big \Vert u_{k}\Big \Vert _{KR}\searrow 
0$ for frames/bases on $ L^{2}(I^{d})$.\ \par 
\, \, Main results are Theorem 3.1, Theorem 3.2, and Theorem 

\bigskip

\noindent{\bf Acknowledgements.} \rm The authors are most grateful 
fo Efim Gluskin of Tel-Aviv University for proposing the scheme of 
proof for Part (2) of Theorem 6.1 and realizing it for the dimension 
$ d=1$; to our regrets, Efim declined our invitation to cosign the 
paper. We also grateful to Sergei Kisliakov of the Steklov Institute, 
St. Petersburg (Russia) for a valuable email exchange on Orlicz space 
interpolation, and to Vasily Vasyunin for carefully reading the manuscript.
\ \par 
\ \par 
  \section{Definitions and comments}
\ \par 
 In order to simplify the statements, we 
\it always assume that our sequences $ (u_{k})_{k\geq 1}$ (frames, 
bases, etc) lay in an one codimensional subspace\ \par 
\ \par 
  \centerline{ $ L^{2}_{0}(\Omega ,m)=\{f\in L^{2}_{{\Bbb R}}(\Omega 
,m):\int _{\Omega }fdm=0\}$\rm .}
\ \par 
\rm The most of results below are still true for all \it Bessel sequences 
$ u=(u_{k})_{k\geq 1}$ \rm in $ L^{2}_{0}$, i.e. the sequences 
with
\ \par 
 $$
  \sum _{k}\Big \vert (f,u_{k})\Big \vert ^{2}\leq 
B(u)^{2}\Vert f\Vert _{2}^{2},\forall f\in L^{2}_{0}
$$
where $ B(u)>0$ stands for the best possible constant in such inequality. 
Recall also that \it a frame \rm (in $ L^{2}_{0}$) is a sequence having\ \par 
\ \par 
  $$ 
  b\Vert f\Vert _{2}^{2}\leq \displaystyle \sum _{k}\Big \vert 
(f,u_{k})\Big \vert ^{2}\leq  B\Vert f\Vert _{2}^{2},\quad \forall f\in L^{2}_{0},
$$
with some constants $ 0<b,B<\infty $, and \it a 
Riesz basis \rm is (by definition) an isomorphic image of an orthonormal 
basis.\ \par 
We always assume that \it the space $ (\Omega ,\rho )$ 
is compact \rm (unless the contrary explicitly follows from the context) 
and \it the measure $ m$ is finite and continuous (has no point masses)\rm .\ \par 
Below, $ \Vert u\Vert _{KR}$ \it stands for the Kantorovich--Rubinstein \rm (also called \it Wasserstein\rm ) \it norm (KR) \rm of 
a zero mean ($ \int udm=0$) signed measure $ udm$; that norm 
evaluates the work needed to transport the positive mass $ u^{+}dm$ 
into the negative one $ u^{-}dm$. In fact, the KR distance $ d(u_{k}^{+}dx,u_{k}^{-}dx)$ 
between measures $ u_{k}^{\pm }dx$ (first invented by L.~Kantorovich 
as early as in 1942, see \cite{K1942}) is a partial case of a more general setting. Namely, 
given nonnegative measures $ \mu ,\nu $ on $ \Omega $ of an equal total 
mass, $ \mu (\Omega )=\nu (\Omega )$, the $ KR$-distance $ d(\mu 
,\nu )$ is defined as the optimal "transfer plan" of the mass distribution 
$ \mu $ to the mass distribution $ \nu $:
\ \par 
$$
 d(\mu ,\nu )=\inf\Big \{\displaystyle \int _{\Omega 
\times \Omega }\rho (x,y)d\psi (x,y):\psi \in \Psi (\mu ,\nu )\Big \},
$$
where the family $ \Psi (\mu ,\nu )$ consists of all "admissible 
transfer plans" $ \psi $, i.e. nonnegative measures on $ \Omega \times \Omega 
$ satisfying the balance (marginal) conditions $ \psi (\Omega \times \sigma 
)-\psi (\sigma \times \Omega )=(\mu -\nu )(\sigma )$ 
for every $ \sigma \subset \Omega $ (the value $ \psi (\sigma \times 
\sigma ')$ has the meaning of how many mass is supposed to transfer 
from $ \sigma $ to $ \sigma '$). The $ KR$-norm of a real (signed) 
measure $ \mu =$ $ \mu _{+}-\mu _{-}$, $ \mu (\Omega )=$ $ 0$, is defined 
as\ \par
\ \par 
  \centerline{ $ \Vert \mu \Vert _{KR}=d(\mu _{+},\mu _{-})$.}
\ \par 
\rm It is shown in Kantorovich--Rubinstein theory (see, for example  \cite{KR1957} or
\cite{KA1977}, Ch.VIII, {\S}4) that the $ KR$-norm of a real (signed) measure 
$ \mu ,\mu (\Omega )=0$, is the dual norm of the Lipschitz 
space\ \par 
\ \par 
  \centerline{ $ \Lip:=\Lip(\Omega )=\{f:\Omega \longrightarrow {\Bbb R}:
\vert f(x)-f(y)\vert \leq c\rho (x,y)\}$}
\ \par 
\rm modulo the constants, where the least possible constant 
$ c$ defines the norm $ \Lip(f)$. Namely,\ \par 
\ \par 
 $$
 \Vert\mu\Vert_{KR}=d(\mu _{+},\mu 
_{-})=\sup\Big \{\displaystyle \int _{I}fd\mu :\Lip(f)\leq 
1\Big \},
$$
where, in fact, it suffices to test only functions $ f\in \text{lip}$, 
$ \text{lip}:=\{f\in \Lip:\vert f(x)-f(y)\vert =o(\rho 
(x,y))\,\,\text{as}\,\,\rho (x,y)\longrightarrow 0\}$. Of course, 
one can extend the above definition to an arbitrary real valued measure 
$ \mu $ setting $ \Vert \mu \Vert =\Vert \mu -\mu (\Omega )\Vert 
_{KR}+\vert \mu (\Omega )\vert $. It makes possible to apply 
our results to $ L^{2}_{{\Bbb R}}$ spaces instead of $ L^{2}_{{\Bbb R},0}$ 
(using that in the case of Bessel sequences, the sequence $ \int _{\Omega }u_{k}dm=
(1,u_{k})$ is in $ l^{2}$). The $ KR$-norm and its variations (with 
various cost function $ h(x,y)$ instead of the distance $ \rho (x,y)$) 
are largely used in the Monge/Kantorovich transportation problems, 
in ergodic theory, etc. We refer to \cite{KA1977} for a basic exposition 
and references, and to \cite{BK2012}, \cite{BKP2017} for extensive and very useful 
surveys of the actual state of the fields.\ \par 
It is clear from the above definitions that, for measuring 
the sign intermixing of $ u_{k}dm$ for a Bessel sequence $ (u_{k})\subset 
L^{2}_{0}$, one can employ certain \it size characteristics 
\rm of the following compact subset of $ L^{2}(\Omega ,m)$,\ \par 
\ \par 
 $$
  \Lip_{1}=\Big \{f:\Omega \longrightarrow {\Bbb R}:
\Big \vert f(x)-f(y)\Big \vert \leq \rho (x,y),f(x_{0})=0\Big \},
$$
\ \par 
\rm where $ x_{0}\in \Omega $ stands for a fixed point of $ \Omega $ 
(it will be easily seen that the choice of $ x_{0}$ does not matter). 
Below, we do that making use of the known Bernstein width numbers $ 
b_{n}(\Lip_{1})$, or - in the case when there exists a linear Hilbert 
space operator $ T$ for which $ \Lip_{1}$ is the range of the unit ball 
- simply the singular numbers $ s_{n}(T)$.\ \par 
Namely, S.Bernstein $ n$-widths $ b_{n}(A,X)$ of a 
(compact) subset $ A\subset X$ (convex, closed and centrally 
symmetric) of a Banach space $ X$ are defined as follows (see \cite{Pi1985}):\ \par 
\ \par 
$$
 b_{n}(A,X)=\sup_{X_{n+1}}\sup\Big \{\lambda :
\lambda B(X_{n+1})\subset A,\lambda \geq 0\Big \},
$$
 where $ X_{n+1}$ runs over all linear subspaces in $ X$ of $ dimX_{n+1}=
n+1$, and $ B(X_{n+1})$ stands for the closed unit ball of $ X_{n+1}$. 
A subspace $ X_{n+1}(A)$ where $ \sup_{X_{n+1}}$ is attained, is called 
optimal; it does not need to be unique (in general). In the case of 
a Hilbert space $ H$ (as everywhere below), if $ A$ is the image of 
the unit ball with respect to a linear (compact) operator $ T$, $ A=
TB(H)$, we have $ b_{n}(A,H)=s_{n}(T)$, where $ s_{k}(T)\searrow 0$ 
($ k=0,1,...$) stands for the $ k$-th singular number of $ T$; 
optimal subspaces $ H_{n+1}(T)$ are simply the linear hulls of $ y_{0},...,y_{n}$ 
from the Schmidt decomposition of $ T$,\ \par 
\ \par 
 $$
  T=\displaystyle \sum _{k\geq 0}s_{k}(T)(\cdot ,x_{k})y_{k},
$$
$ (x_{k})$ \rm and $ (y_{k})$ being orthonormal sequences in $ H$.\ \par 
\ \par 
  \section{ Statements}
\ \par 
 Recall that $ (\Omega ,\rho )$ stands for a compact 
metric space (unless the other is claimed explicitly), and $ m$ is a 
finite Borel measure on $ \Omega $ \it having no point masses \rm (for 
convenience normalized to $ 1$).\ \par 
Lemma 1 below shows what kind of the intermixing of 
signs we have for free, for every Bessel sequence $ (u_{k})$. Lemma 
2 shows that in no cases, one can have an intermixing better than $ 
l^{2}$ smallness of $ \Vert u_{k}\Vert _{KR}$. All intermediate cases 
can occur, following the widths properties of the compact $ \Lip_{1}\subset 
L^{2}(\Omega ,m)$, see Theorems 3.1,3.2 and the comments below.\ \par 
\ \par 
\bf Lemma 1. \it For every Bessel sequence $ (u_{k})_{k\geq 
1}$ in $ L^{2}_{{\Bbb R}}(\Omega ,m)$, we have\ \par 
\ \par 
  \centerline{ $ \lim_{k}\Vert u_{k}\Vert _{KR}=0$.}
\ \par 
\bf Lemma 2. \it For every compact measure triple $ 
(\Omega ,\rho ,m)$ (with the above conditions) and every sequence $ 
(\epsilon _{k})_{k\geq 1}$, $ \epsilon _{k}\geq 0$, such that $ \sum _{k}\epsilon 
_{k}^{2}<\infty $, there exists an orthonormal sequence $ (u_{k})_{k\geq 
1}$ in $ L^{2}_{{\Bbb R}}(\Omega ,m)$ satisfying\ \par 
\ \par 
  \centerline{ $\Vert u_{k}\Vert _{KR}\geq c\epsilon _{k},
k=1,2,...$ \rm ($ c>0$)\it .}
\ \par 
\it In particular, there exists an orthonormal sequence $ (u_{k})_{k\geq 
1}$ in $ L^{2}_{{\Bbb R}}(\Omega ,m)$ such that\ \par 
\ \par 
  \centerline{ $ \displaystyle \sum _{k}\Vert u_{k}\Vert ^{2-\epsilon 
}_{KR}=\infty ,\forall \epsilon >0$.}
\ \par 
\bf Lemma 3. \it For every sequence $ (\epsilon _{k})_{k\geq 
1}$, $ \epsilon _{k}>0$, with $ \lim_{k}\epsilon _{k}=0$, there 
exists a compact measure triple $ (\Omega ,\rho ,m)$ (with the above 
conditions) and an orthonormal sequence $ (u_{k})_{k\geq 1}$ in $ L^{2}_{{\Bbb R}}(\Omega 
,m)$ such that\ \par 
\ \par 
  \centerline{ $\Vert u_{k}\Vert _{KR}=c\epsilon _{k},
k=1,2,...$ \rm ($ {\frac{1}{2\sqrt{2} }} \leq c\leq {\frac{\displaystyle 2\sqrt{\displaystyle 
2} }{\displaystyle \pi }} $).}
\ \par 
\ \par 
\bf Theorem 3.1. \it (1) Given a Bessel sequence $ (u_{k})_{k\geq 
1}$ in $ L^{2}_{{\Bbb R}}(I,dx)$, $ I=(0,1)$\rm , \it we have\ \par 
\ \par 
  \centerline{ $ \displaystyle \sum _{k}\Vert u_{k}\Vert _{KR}^{2}<
\infty .$}
\ \par 
\it (2) Given a Bessel sequence $ (u_{k})_{k\geq 1}$ in $ L^{2}_{{\Bbb R}}(I^{d},dx)$, 
$ d=2,3,...$, we have\ \par 
\ \par 
  \centerline{ $ \displaystyle \sum _{k}\Vert u_{k}\Vert _{KR}^{d+\epsilon 
}<\infty $\rm , $ \forall \epsilon >0$\sc .}
\ \par 
\it (3) For the $ Sin$ orthonormal sequence $ (u_{n})_{n\in 2{\Bbb N}^{d}}$ 
in $ L^{2}_{{\Bbb R}}(I^{d},dx)$,\ \par 
\ \par 
  \centerline{ $ u_{n}(x)=2^{d/2}Sin(\pi n_{1}x_{1})Sin(\pi n_{2}x_{2})...Sin(\pi 
n_{d}x_{d})$ ($ n=(n_{1},...,n_{d})\in (2{\Bbb N})^{d}$),}
\ \par 
\it we have\ \par 
\ \par 
  \centerline{ $ \displaystyle \sum _{n}\Vert u_{n}\Vert _{KR}^{d}=
\infty $.}
\ \par 
\bf Remark. \rm For a generic Bessel sequence (or, an orthonormal sequence), 
the $ l^{2}$-convergence property (1) is a best possible result (see 
Lemma 2). However, for certain specific sequences, \it (1) can be much 
sharpen\rm . For example, let $ u\in L^{2}_{{\Bbb R},0}({\Bbb T})$ 
and\ \par 
\ \par 
  \centerline{ $ u_{n}(\zeta)=u(\zeta^{n})$, $ n=1,2,...$}
\ \par 
\rm Then, as it easy to see,\ \par 
\ \par 
  \centerline{ $ \Vert u_{n}\Vert _{KR}\leq {\frac{1}{n}} \Vert 
u\Vert _{KR}$ }
\ \par 
\rm (in fact, there is an equality), and so $ \displaystyle \sum _{n}\Vert 
u_{n}\Vert _{KR}^{1+\epsilon }<$ $ \infty $ ($ \forall \epsilon >0$). 
Such a dilated sequence $ (u_{n})_{n}$ is Bessel if, and only if, the 
\it Bohr transform \rm of $ u$, $ Bu(\zeta)=\sum _{n}\hat u(n)\zeta^{\alpha 
(n)}$, $ \zeta^{\alpha }=\zeta_{1}^{\alpha _{1}}\zeta^{\alpha _{2}}...$ 
($ n=2^{\alpha _{1}}3^{\alpha _{2}}...$ stands for for Euclid 
prime representation of $ n\in {\Bbb N}$) is bounded on the multitorus 
$ \zeta=(\zeta_{1},\zeta_{2},...)\in {\Bbb T}^{\infty }$, see 
for instance \cite{Ni2017}.\ \par 
\ \par 
In fact, Theorem 3.1, is an immediate corollary 
of the next Theorem 3.2. We extend the property $ (\Vert u_{k}\Vert _{KR})\in 
l^{2}$ to any "one dimensional smooth manifolds", see Proposition 5.1 
for the exact statement. Lemma 2 shows that this condition describe 
the fastest decrease of the $ KR$-norms for a generic Bessel sequence. 
On the spaces $ \Omega ,\rho $ of "higher dimensions" the property 
fails.\ \par 
In Theorem 3.2, we develop the approach mentioned at 
the end of Section 2: we compare the compact set $ \Lip_{1}$ with the 
$ T$-range $ T(B(L^{2}))$of the unit ball for an appropriate compact 
operator $ T$. For a direct comparison $\Vert u_{n}\Vert _{KR}$ 
with Bernstein numbers $ b_{n}(\Lip_{1})$ see Section 5 below.\ \par 
\ \par 
\bf Theorem 3.2. \it Let $ T$ be compact linear operator\ \par 
\ \par 
  \centerline{ $ T:L^{2}_{{\Bbb R}}(\Omega ,m)\longrightarrow L^{2}_{{\Bbb R}}(\Omega 
,m)$,}
\ \par 
\it and $ \varphi :[0,\infty )\longrightarrow [0,\infty )$ 
be a continuous increasing function on $ [0,\infty )$ whose inverse 
$ \varphi ^{-1}$ satisfies 
\ \par 
$$
 \varphi ^{-1}(x)=x^{1/2}r(1/x^{-1/2})\quad \forall 
x>0
$$
 with a concave (or, pseudo-concave) function $ x\longmapsto r(x)$ 
on $ (0,\infty )$.\ \par 
(1) If $ \Lip_{1}\subset $ $ T(B(L^{2}_{{\Bbb R}}(\Omega ,m)))$ and 
$ \sum _{k}\varphi (s_{k}(T))<\infty $, then, for every Bessel 
sequence $ (u_{k})\subset L^{2}_{{\Bbb R}}(\Omega ,m)$,\ \par 
\ \par 
  \centerline{ $ \sum _{k\geq 1}\varphi (a\Vert u_{k}\Vert _{KR})<
\infty $ (for a suitable $ a>0$).}
\ \par 
\it (2) If $ \Lip_{1}\supset $ $ T(B(L^{2}_{{\Bbb R}}(\Omega ,m)))$, 
then there exists an orthonormal sequence $ (u_{k})_{k\geq 0}\subset $ 
$ L^{2}_{{\Bbb R}}(\Omega ,m)$, such that\ \par 
\ \par 
  \centerline{ $ \Vert u_{k}\Vert _{KR}\geq s_{k}(T)$, $ k=
0,1,...$}
\ \par 
\it In particular (in order to compare with (1)), $ \sum _{k}h(\Vert u_{k}\Vert 
_{KR})=$ $ \infty $ for every $ h$ for which $ \sum _{k}h(s_{k}(T))=
\infty $.\ \par 
\ \par 
\bf Remark. \rm See Section 5.III below for a version of Theorem 3.2, 
point (2), employing the Bernstein widths $ b_{n}(\Lip_{1})$ instead 
of $ s_{n}(T)$ ($ T$ does not need to exist for the compact set $ \Lip_{1}$).\ \par 
\ \par 
{\bf Corollary.} {\it Let $ \Lip_{1}=T(B(L^{2}_{{\Bbb R}}(\Omega 
,m)))$ and $ p(T):=\inf\{\alpha :\sum _{k}s_{k}(T)^{\alpha }<$ 
$ \infty \}$.\ \par 
(1) If $ p(T)<2$, then $ \sum _{k}\Vert u_{k}\Vert ^{2}<$ $ \infty $, 
for every Bessel sequence $ (u_{k})\subset $ $ L^{2}_{{\Bbb R}}(\Omega ,m)$. 
On the other hand, there exists $ T$ with $ p(T)=1$ and an 
orthonormal sequence such that $ \sum _{k}\Vert u_{k}\Vert _{KR}^{2-\epsilon 
}=$ $ \infty $ ($ \forall \epsilon >0$) (see Lemma 2 above)\ \par 
(2) If $ \sum _{k}s_{k}(T)^{p}<\infty $, $ p\geq 2$, then $ 
\sum _{k}\Vert u_{k}\Vert _{KR}^{p}<\infty $ for every Bessel 
sequence $ (u_{k})\subset $ $ L^{2}_{{\Bbb R}}(\Omega ,m)$.}\ \par 
\ \par 
{\bf Remark.} As we will see, Theorem 3.1, in fact, is a consequence 
of the last Corollary. Some concrete examples to Theorem 3.2 are presented 
below, in Section 5.\ \par 
\ \par 
  \section{ Proofs}
\label{proofs}
\bf I. Proof of Lemma 1. 
%\begin{proof}
\rm Since $ (u_{k})_{k\geq 1}$ is a Bessel 
sequence, it tends weakly to zero: $ (u_{k},f)\longrightarrow 0$ 
as $ k\longrightarrow \infty $, for every $ f\in L^{2}_{{\Bbb R}}(\Omega 
,m)$. On a (pre)compact set $ f\in \Lip_{1}$, the limit is uniform:\ \par  
 $$
  \lim_{k}\Vert u_{k}\Vert _{KR}= \lim_{k}\sup\Big \{\displaystyle \int 
_{\Omega }u_{k}fd\mu : f\in \Lip_{1}\Big \}=\,0\,.
$$
\ \par
%\end{proof} 
\bf II. Proof of Lemma 2.
%\begin{proof}
 \rm The Borel measure $ m$ being continuous 
satisfies the Menger property: the values $ mE$, $ E\subset \Omega $ 
fill in an interval $ [0,m(\Omega )]$; if $ m$ is normalized - the 
interval $ [0,1]$ (see \cite{Ha1950}, {\S}41 (with many retrospective references, 
the oldest one is to K.Menger, 1928), and for a complete and short 
proof \cite{DN2011}, Prop. A1, p.645). Below, we use that property many 
times.\ \par 
Let $ E_{i}\subset \,\Omega $ be disjoint Borel 
sets, $ E_{1}\bigcap E_{2}=\,\emptyset $, $ mE_{i}=\,1/2$, 
and further, $ K_{i}\subset \,E_{i}$ be compacts such that $ mK_{i}=\,
1/3$ ($ i=\,1,2$). Denote $ \delta =\,dist(K_{1},K_{2})>\,
0$, and set\ \par 
\ \par 
  \centerline{ $ f(x)=(1-{\frac{2}{\delta }} dist(x,K_{1}))^{+
}-(1-{\frac{2}{\delta }} dist(x,K_{2}))^{+}$, $ x\in \Omega $.}
\ \par 
\rm Then, $ f\in \Lip(\Omega ,\rho )$, $ \Lip(f)\leq 2/\delta $ 
and $ f(x)=1$ for $ x\in K_{1}$, $ f(x)=-1$ for $ x\in 
K_{2}$.\ \par 
Now, using the Menger property, one can find two sequences 
$ (\Delta _{k}^{1})$, $ (\Delta _{k}^{2})$, $ k=1,2,...$, of pairwise 
disjoint sets such that $ \Delta ^{i}_{k}\subset K_{i}$, $ \Delta 
^{i}_{k}\bigcap \Delta ^{i}_{j}=\emptyset $ ($ i=1,2$, 
$ k\not= j$), and $ m\Delta ^{1}_{k}=m\Delta ^{2}_{k}=a^{2}\epsilon 
_{k}^{2}$, where $ a>0$ is chosen in such a way that $ a^{2}\sum _{k\geq 
1}\epsilon _{k}^{2}\leq 1/3$. Setting
\ \par 
 $$
  u_{k}=c_{k}(\chi _{\Delta ^{1}_{k}}-\chi 
_{\Delta ^{2}_{k}}),\quad k=1,2,...,
$$
 with $ \Vert u_{k}\Vert _{2}^{2}=2c_{k}^{2}m\Delta ^{1}_{k}=
1$, we obtain an orthonormal sequence $ (u_{k})\subset L^{2}(\Omega 
,m)$ such that\ \par 
\ \par 
$\Vert u_{k}\Vert _{KR}\geq \displaystyle \int _{\Omega 
}u_{k}({\frac{\displaystyle \delta }{\displaystyle 2}} f)dm={\frac{\displaystyle 
\delta }{\displaystyle 2}} 2c_{k}m\Delta ^{1}_{k}={\frac{\displaystyle \delta 
}{\displaystyle \sqrt{\displaystyle 2} }} \sqrt{\displaystyle m\Delta ^{1}_{k}} 
={\frac{\displaystyle \delta a}{\displaystyle \sqrt{\displaystyle 2} 
}} \epsilon _{k}$. \par 
%\end{proof}
\bf III. Proof of Lemma 3.
%\begin{proof}
 \rm Let $ \Omega ={\Bbb T}^{\infty }$, 
the infinite topological product of compact abelian groups $ {\Bbb T}\times 
{\Bbb T}\times ...$, endowed with its normalized Haar measure $ m_{\infty }=
m\times m\times ...$. The product topology on $ \Omega $ is metrizable 
by a variety of metrics, we choose $ \rho =\rho _{\epsilon }$, 
$ \epsilon =(\epsilon _{k})_{k\geq 1}$ defined by\ \par 
\ \par 
  \centerline{ $ \rho _{\epsilon }(\zeta,\zeta')=max_{k\geq 1}\epsilon 
_{k}\vert \zeta_{k}-\zeta'_{k}\vert $, $ \zeta',\zeta=(\zeta_{k})_{k\geq 
1}\in {\Bbb T}^{\infty }$.}
\ \par 
\rm Setting
\ \par 
 $$
  u_{k}(\zeta)=\sqrt{2} Re(\zeta_{k}), \quad \zeta\in 
{\Bbb T}^{\infty },
$$
we define an orthonormal sequence in $ L^{2}({\Bbb T}^{\infty },m_{\infty 
})$ with $ \vert u_{k}(\zeta)-u_{k}(\zeta')\vert \leq {\frac{\sqrt{2} 
}{\epsilon _{k}}} \rho (\zeta,\zeta')$, and so $ \Lip(u_{k})\leq \sqrt{2} 
/\epsilon _{k}$.\ \par 
Further, we need the following notation: let 
$ f\in \Lip_{1}({\Bbb T}^{\infty })$, $ f(\zeta)=$ $ f(\zeta_{k},\overline{\zeta})$ 
where $ \zeta=$ $ (\zeta_{k},\overline{\zeta})\in {\Bbb T}^{\infty }=$ 
$ {\Bbb T}\times {\Bbb T}^{\infty }$, $ \overline{\zeta}$ consists 
of variables different from $ \zeta_{k}$, and\ \par 
\ \par 
$$
 \overline{u}_{k}(\zeta_{k})= \sqrt{2} Re(\zeta_{k}), \quad \zeta_{k}\in {\Bbb T},
$$
(in fact, this is one and the same function $ e^{i\theta }\longmapsto 
\sqrt{2} Cos(\theta )$ for every $ k$). Finally, we set $ \overline{f}(\zeta_{k}):=
\displaystyle \int _{{\Bbb T}^{\infty }}f(\zeta_{k},\overline{\zeta})dm_{\infty 
}(\overline{\zeta})$ and observe that $ \Lip(\overline{f})\leq \epsilon 
_{k}$:\ \par 
\ \par 
  \centerline{ $ \Big \vert \overline{f}(\zeta_{k})-\overline{f}(\zeta'_{k})\Big \vert 
\leq $ $ \displaystyle \int _{{\Bbb T}^{\infty }}\Big \vert f(\zeta_{k},\overline{\zeta})-f(\zeta'_{k},\overline{\zeta})\Big \vert 
dm_{\infty }(\overline{\zeta})\leq $}
\ \par 
  \centerline{ $ \leq \displaystyle \int _{{\Bbb T}^{\infty }}\epsilon 
_{k}\Big \vert \zeta_{k}-\zeta'_{k}\Big \vert dm_{\infty }(\overline{\zeta})
=\epsilon _{k}\Big \vert \zeta_{k}-\zeta'_{k}\Big \vert $\rm .}
\ \par 
\rm Now,\ \par 
  \centerline{ $ \displaystyle \int _{{\Bbb T}^{\infty }}u_{k}(\zeta)f(\zeta_{k},\overline{\zeta})dm_{\infty 
}(\zeta)=\displaystyle \int _{{\Bbb T}}\overline{u}_{k}(\zeta_{k})\displaystyle \int 
_{{\Bbb T}^{\infty }}f(\zeta_{k},\overline{\zeta})dm_{\infty }(\overline{\zeta})dm(\zeta_{k})=$}
\ \par 
  \centerline{ $ =\displaystyle \int _{{\Bbb T}}\overline{u}_{k}(\zeta_{k})\overline{f}(\zeta_{k})dm(\zeta_{k})\leq 
\epsilon _{k}\Vert \overline{u}_{k}\Vert _{KR({\Bbb T})}$\rm ,}
\ \par 
\rm and hence $\Vert u_{k}\Vert _{KR({\Bbb T}^{\infty }  )}\leq 
\epsilon _{k}\Vert \overline{u}_{k}\Vert _{KR({\Bbb T})}$.\ \par 
Conversely, if $ h\in \Lip_{1}({\Bbb T})$ and $ \underline{h}(\zeta):=
h(\zeta_{k})$ for $ \zeta\in {\Bbb T}^{\infty }$, then $ \Big \vert \underline{h}(\zeta_{k})-\underline{h}(\zeta'_{k})\Big \vert 
\leq {\frac{\displaystyle 1}{\displaystyle \epsilon _{k}}} \rho 
(\zeta,\zeta')$, and so\ \par 
\ \par 
  \centerline{ $ \displaystyle \int _{{\Bbb T}}\overline{u}_{k}hdm(\zeta_{k})=$ 
$ \displaystyle \int _{{\Bbb T}^{\infty }}dm_{\infty }(\overline{\zeta})\displaystyle \int 
_{{\Bbb T}}\overline{u}_{k}(\zeta_{k})h(\zeta_{k})dm(\zeta_{k})=\displaystyle \int 
_{{\Bbb T}^{\infty }}u_{k}(\zeta)\underline{h}(\zeta)dm_{\infty }(\zeta)\leq 
$}
\ \par 
  \centerline{ $ \leq {\frac{\displaystyle 1}{\displaystyle \epsilon 
_{k}}}\Vert u_{k}\Vert _{KR({\Bbb T}^{\infty }  )}$\rm ,}
\ \par 
\rm which entails $\Vert\overline{u}_{k}\Vert _{KR({\Bbb T})}\leq 
{\frac{\displaystyle 1}{\displaystyle \epsilon _{k}}}\Vert u_{k}\Vert 
_{KR({\Bbb T}^{\infty }  )}$. Finally, $\Vert u_{k}\Vert _{KR({\Bbb T}^{\infty 
}  )}=\epsilon _{k}\Vert \overline{u}_{k}\Vert _{KR({\Bbb T})}$.\ \par 
Moreover, since $ \Lip(\overline{u}_{k})\leq \sqrt{2} 
$,\ \par 
\ \par 
  \centerline{ $ {\frac{\displaystyle 1}{\displaystyle 2\sqrt{\displaystyle 2} 
}} =\displaystyle \int _{{\Bbb T}}\overline{u}_{k}(\overline{u}_{k}/\sqrt{\displaystyle 
2} )dm(\zeta_{k})\leq\Vert\overline{u}_{k}\Vert _{KR({\Bbb T})}\leq 
\Vert \overline{u}_{k}\Vert _{L^{1}  ({\Bbb T})}=
{\frac{\displaystyle 2\sqrt{\displaystyle 2} }{\displaystyle \pi }} $. 
}
\ \par 
\bf Remark. \rm For the same space $ L^{2}{\rm (}{\Bbb T}^{\infty }{\rm ,m}_{\infty 
}{\rm )}$, but with a \it non-compact \rm (bounded) metric $ \rho {\rm (}\zeta{\rm ,}\zeta{\rm ')=}$ 
$ \sup_{k\geq 1}\vert \zeta_{k}{\rm -}\zeta{\rm '}_{k}\vert $, we have 
$ \Vert {\rm u}_{{\rm k}}\Vert _{{\rm KR}}\geq $ $ {\rm 1}$ for $ {\rm u}_{{\rm k}}{\rm (}\zeta{\rm )=}$ 
$ {\rm Sin}\pi {\rm x}_{{\rm k}}$, $ \zeta{\rm =}$ $ {\rm (e}^{{\rm ix}_{{\rm 1}}}{\rm ,e}^{{\rm ix}_{{\rm 2}}}{\rm ,...,e}^{{\rm ix}_{{\rm k}}}{\rm ,...)}\in 
{\Bbb T}^{\infty }$, so that $ (\Vert {\rm u}_{k}\Vert _{{\rm KR}})_{k\geq 
1}$ does not tend to zero.\ \par 
%\end{proof}
\bf IV. Proof of Theorem 3.1.
%\begin{proof}
 \it (1) \rm Since $ u_{k}\in L^{2}_{{\Bbb R},0}(I,dx)$, 
$ \int _{I}u_{k}dx=0$. Taking a smooth function $ f$ with $ 
\Lip(f)\leq 1$ (which are dense in the unit ball of $ lip$) and $ v_{k}(x)=
Ju_{k}(x):=\int _{0}^{x}u_{k}dx$, we get $ v_{k}(0)=v_{k}(1)=
0$, and hence\ \par 
\ \par 
  \centerline{ $ \displaystyle \int _{I}fu_{k}dx=(fv_{k})_{0}^{1}-
\displaystyle \int _{I}v_{k}f'dx=-\displaystyle \int _{I}v_{k}f'dx$.}
\ \par 
\rm Making $ \sup$ over all $ f$ with $ \vert f'\vert \leq 1$, we obtain 
$\Vert u_{k}\Vert _{KR}=\Vert v_{k}\Vert _{L^{1}}$. 
But the mapping\ \par 
 $$
  J:L^{2}(I)\longrightarrow L^{2}(I)
$$
is a Hilbert-Schmidt operator, and hence $ \displaystyle \sum _{k}\Vert 
Ju_{k}\Vert ^{2}_{L^{2}}<\infty $, and so $ \displaystyle  \sum _{k}\Vert u_{k}\Vert ^{2}_{KR}= \sum _{k}\Vert 
Ju_{k}\Vert ^{2}_{L^{1}}<
\infty $.\ \par 
The penultimate inequality is obvious if $ (u_{k})$ 
is an orthonormal (or only Riesz) sequence, but is still true for every 
Bessel sequence $ (u_{k})_{k\geq 1}$. Indeed, taking an auxiliary orthonormal 
basis $ (e_{j})_{j\geq 1}$ in $ L^{2}_{{\Bbb R}}(I,dx)$, we can write\ \par 
\ \par 
  \centerline{ $ \displaystyle \sum _{k}\Vert Ju_{k}\Vert ^{2}_{L^{2}}=
\displaystyle \sum _{k}\displaystyle \sum _{j}\Big \vert (Ju_{k},e_{j})\Big \vert 
^{2}=\displaystyle \sum _{j}\displaystyle \sum _{k}\Big \vert (u_{k},J^{*}e_{j})\Big \vert 
^{2}\leq $}
\ \par 
  \centerline{ $ \leq \displaystyle \sum _{j}const\cdot\Vert J^{*}e_{j}\Vert 
^{2}<\infty $\rm ,}
\ \par 
\rm since the adjoint $ J^{*}$ is a Hilbert-Schmidt operator. \ \par 
\it (2) \rm This is a $ d$-dimensional version of the previous reasoning. 
Anew, we use the dual formula for the KR norm,\ \par 
\ \par 
$$
 \Vert u_{k}\Vert _{KR}=\sup\{\int _{I^{d}}fu_{k}dx:
f\in C^{\infty },\Lip(f)\leq 1,\int fdx=0\},
$$
the last requirement does not matter since $ \Lip(f)=\Lip(f+const)$. 
Notice that for $ f\in C^{\infty }(I^{d})$, $ \Lip(f)\leq $ $ 1$ $ \Leftrightarrow 
\vert \nabla f(x)\vert \leq 1$ ($ x\in I^{d}$), where 
$ \nabla f$ stands for the gradient vector $ \nabla f=({\frac{\partial 
f}{\partial x_{j}}} )_{1\leq j\leq d}$ . Now, define a linear mapping 
on the set $ {\cP}_{0}$ of vector valued trigonometric polynomials 
of the form $ \sum _{n\in {\Bbb Z}^{d}}c_{n}\nabla e^{i(n,\cdot )}\in L^{2}(I^{d},{\Bbb C}^{d})$ 
with the zero mean ($ c_{0}=0$) by the formula\ \par 
\ \par 
  \centerline{ $ A(\nabla e^{i(n,x)})=\Big \vert n\Big \vert e^{i(n,x)}$, 
$ n\in {\Bbb Z}^{d}\backslash \{0\}$.}
\ \par 
\rm It is clear that $ A$ extends to a unitary operator\ \par 
\ \par 
  \centerline{ $ A:clos_{L^{2}  (I^{d}  ,{\Bbb C}^{d}  )}(\nabla {\cP}_{0})\longrightarrow 
L_{0}^{2}(I^{d})$.}
\ \par 
\rm Further, let $ M:L_{0}^{2}(I^{d})\longrightarrow L_{0}^{2}(I^{d})$ 
be a (bounded) multiplier,\ \par 
\ \par 
  \centerline{ $ M(e^{i(n,x)})={\frac{\displaystyle 1}{\displaystyle \vert 
n\vert }} e^{i(n,x)}$, $ n\in {\Bbb Z}^{d}\backslash \{0\}$\it ,}
\ \par 
\rm and finally, $ T(\nabla f)=f$, $ f\in C_{0}^{\infty }(I^{d})$. 
Then,\ \par 
\ \par 
  \centerline{ $ \displaystyle \int _{I^{d}}fu_{k}dx=\displaystyle \int _{I^{d}}(T(\nabla 
f))u_{k}dx=\displaystyle \int _{I^{d}}\nabla f\cdot (T^{*}u_{k})dx$,}
\ \par 
$ T^{*}$ \rm being the adjoint between $ L^{2}$ Hilbert spaces. It 
follows\ \par 
\ \par 
  \centerline{ $\Vert u_{k}\Vert _{KR}\leq \sup\Big \{\displaystyle \int 
_{I^{d}}\nabla f(T^{*}u_{k})dx:\Big \vert \nabla f\Big \vert \leq 
1\Big \}\leq\Vert T^{*}u_{k}\Vert _{L^{1}  (I^{d}  ,{\Bbb C}^{d}  )}\leq 
$}
\ \par 
  \centerline{ $ \leq\Vert T^{*}u_{k}\Vert _{L^{2}  (I^{d}  ,{\Bbb C}^{d}  )}$\rm .}
\ \par 
\rm Moreover, $ T=MA$, where $ A$ is unitary (between the corresponding 
spaces) and $ M$ in a Schatten-von Neumann class $ {\cS}_{p}$ for 
every $ p$, $ p>d$ (since $ M$ is diagonal and $ \displaystyle \sum _{n\in 
{\Bbb Z}^{d}  \backslash \{0\}}{\frac{\displaystyle 1}{\displaystyle \vert 
n\vert ^{p}}} <\infty \Leftrightarrow p>d$). 
Using the dual definition of the Bessel sequence as $ \Vert \sum a_{k}u_{k}\Vert 
^{2}\leq $ $ c\,(\sum a_{k}^{2})$ for every real finite sequence $ (a_{k})$, 
we can write $ (u_{k})$ as the image $ u_{k}=Be_{k}$ of an 
orthonormal sequence $ (e_{k})$ under a linear bounded map $ B$. This 
gives\ \par 
\ \par 
  \centerline{ $\Vert u_{k}\Vert _{KR}\leq\Vert T^{*}Be_{k}\Vert 
_{L^{2}}$.}
\ \par 
\rm For every $ p>d$, this implies $ \displaystyle \sum _{k}\Vert u_{k}\Vert 
^{p}_{KR}\leq $ $ \displaystyle \sum _{k}\Vert T^{*}Be_{k}\Vert ^{p}_{L^{2}}<
\infty $ since $ T^{*}B\in {\cS}_{p}$ and $ d\geq 2$ (see Remark 
below). 

\bf Remark. \rm For the last property, see for example \cite{GoKr1965}. 
Here is a simple explanation: given a linear bounded operator $ S:H\longrightarrow 
K$ between two Hilbert spaces and an orthonormal sequence $ (e_{k})$ 
in $ H$, define a mapping $ j:S\longrightarrow (Se_{k})$; then, 
$ j$ is bounded as a map $ {\cS}_{2}\longmapsto l^{2}(K)$ and as 
a map $ {\cS}_{\infty }\longmapsto c_{0}(K)$ (compact operators); 
by operator interpolation, $ j:{\cS}_{p}\longmapsto l^{p}(K)$ is 
also bounded for $ 2<p<\infty $. 
\ \par 
For $ 1\leq p\leq 2$, the things go differently: the 
best summation property $ \displaystyle \sum _{k}\Vert Se_{k}\Vert ^{\alpha 
}<$ $ \infty $, which one can generally have for $ S\in {\cS}_{p}$, 
is only for $ \alpha =2$ (look at rank one operators $ S=(\cdot 
,x)y$). This claim explains the strange behavior in exponent from 
$ 2+\epsilon $ for dimension $ 2$ to exactly $ 2$ for dimension $ 1$ 
(and not $ 1+\epsilon $ as one would expect).\ \par 
\ \par 
\it (3) \rm We use anew the duality formula\ \par 
\ \par 
  \centerline{ $\Vert u_{n}\Vert _{KR}=$ $ \sup\Big \{\displaystyle \int _{I^{d}}fu_{n}d\mu 
:$ $ \Lip(f)\leq $ $ 1\Big \}$.}
\ \par 
\rm Taking $ f=u_{n}/\Lip(u_{n})$ we get $\Vert u_{n}\Vert _{KR}\geq 
1/\Lip(u_{n})$ where $ \Lip(u_{n})\leq max\vert \nabla u_{n}(x)\vert 
\leq 2^{d/2}\vert n\vert $, and so \ \par 
\ \par 
  \centerline{ $ \displaystyle \sum _{n}\Vert u_{n}\Vert _{KR}^{d}\geq 
2^{-d^{2}  /2}\displaystyle \sum _{n\in (2{\Bbb N})^{d}}\Big \vert n\Big \vert 
^{-d}=\infty $. }
%\end{proof}

\bf V. Proof of Theorem 3.2. 
%\begin{proof}
\rm Let $ T=\sum _{k\geq 0}s_{k}(T)(\cdot 
,x_{k})y_{k}$ be the Schmidt decomposition of a compact operator $ T$ 
acting on a Hilbert space $ H$, $ s_{k}(T)\searrow 0$ being the singular 
numbers. Let further, $ A:H\longrightarrow H$ be a bounded operator, 
and $ (e_{k})_{k\geq 0}$ an arbitrary (fixed) orthonormal basis. Given 
a sequence $ \alpha =(\alpha _{j})_{j\geq 0}$ of real numbers, 
$ \alpha \in l^{\infty }$, define a bounded operator\ \par 
\ \par 
  \centerline{ $ T_{\alpha }=$ $ \sum _{k\geq 0}\alpha _{k}(\cdot ,x_{k})y_{k}$\it ,}
\ \par 
\rm and then a mapping
\ \par 
$$ 
j:\alpha \longmapsto (T_{\alpha }^{*}Ae_{k})_{k\geq 
0},
$$
a $ H$-vector valued sequence in $ l^{\infty }(H)$.\ \par 
We are using a (partial case of a) J.~Gustavsson--J.~Peetre 
interpolation theorem \cite{GuP1977}for Orlicz spaces. Recall that, in 
the case of sequence spaces, an Orlicz space $ l^{\varphi }$, where 
$ \varphi :{\Bbb R}_{+}\longrightarrow {\Bbb R}_{+}=(0,\infty )$ 
is increasing, continuous, and meets the so-called $ \Delta _{2}$-condition $ \varphi (2x)\leq 
C\varphi (x)$, $ x\in {\Bbb R}_{+}$, is the vector space of 
real sequences $ c=(c_{k})$ satisfying $ \sum _{k}\varphi (a\vert 
c_{k}\vert )<\infty $ for a suitable $ a>0$. Similarly, a vector 
valued Orlicz space consists of sequences $ c=$ $ (c_{k})$, $ c_{k}\in 
H$ having $ \sum _{k}\varphi (a\Vert c_{k}\Vert )<$ $ \infty $ for 
a suitable $ a>0$. We need the Hilbert space valued spaces only. The 
Gustavsson--Peetre interpolation theorem (theorem 9.1 in \cite{GuP1977}) 
implies that if mappings $ j:l^{\infty }\longrightarrow l^{\infty }(H)$ 
and $ j:l^{2}\longrightarrow l^{2}(H)$ are bounded, then
\ \par 
$$
 j:l^{\varphi }\longrightarrow l^{\varphi }(H)
$$
is bounded whenever the measuring function $ \varphi $ satisfies 
the conditions given in Theorem 3.2.\ \par 
\it (1) \rm Now, in the notation and the assumptions of statement (1), 
the Bessel sequence $ (u_{k})$ is of the form $ u_{k}=Ae_{k}$, 
where $ A$ is a bounded operator and $ (e_{k})$ an orthonormal sequence. 
It follows\ \par 
\ \par 
  \centerline{ $\Vert u_{k}\Vert _{KR}=\sup_{f\in \Lip_{1}}\Big \vert 
(Ae_{k},f)\Big \vert \leq \sup_{f\in T(B(L^{2}  ))}\Big \vert (Ae_{k},f)_{L^{2}}\Big \vert 
=\Vert T^{*}Ae_{k}\Vert _{L^{2}}$.}
\ \par 
\rm For every $ \alpha \in l^{2}$, $ T_{\alpha }\in {\cS}_{2}$ 
(Hilbert-Schmidt), and then $ T_{\alpha }^{*}A\in {\cS}_{2}$, and 
hence $ j(\alpha )\in l^{2}(H)$. By Gustavsson--Peetre, $ \alpha \in l^{\varphi 
}\Rightarrow j(\alpha )\in l^{\varphi }(H)$. Applying 
this with $ \alpha =(s_{k}(T))$, we get $ \sum _{k}\varphi (a\Vert 
u_{k}\Vert _{KR})\leq \sum _{k}\varphi (a\Vert T^{*}Ae_{k}\Vert )<$ 
$ \infty $ for a suitable $ a>0$. 
\ \par 
\it (2) \rm In the assumptions of (2), and with the Schmidt decomposition 
$$ 
T= \sum _{k\geq 0}s_{k}(T)(\cdot ,x_{k})y_{k}\,,
$$
 set $ u_{k}=y_{k}$
$ k\geq 0$. Then\ \par 
\ \par 
$\Vert u_{k}\Vert _{KR}=$ $ \sup_{f\in \Lip_{1}}\Big \vert (y_{k},f)\Big \vert 
\geq $ $ \sup_{f\in T(B(L^{2}  ))}\Big \vert (y_{k},f)\Big \vert =$ 
$\Vert T^{*}y_{k}\Vert _{2}=s_{k}(T)$. 
\ \par 
%\end{proof}
  \section{Further examples and comments}
\ \par 
\bf I. Fastest and slowest rates of decreasing $\Vert u_{k}\Vert _{KR}\searrow 
0$. \rm Lemma 2 shows that, the $ KR$-norms of a generic Bessel sequence 
don't have to be smaller than required by the condition $ \displaystyle \sum _{k}\Vert 
u_{k}\Vert ^{2}_{KR}<\infty $.\ \par 
On the other hand, point (1) of Theorem 3.1 gives an 
example of $ (\Omega ,\rho ,dx)$, where every Bessel sequence meets 
that property. 

Now, we extend this 
result to measure spaces over (almost) arbitrary 1-dimensional "smooth 
manifold" of finite length, as follows.\ \par 

\medskip

As to the fastest possible decreasing of $ \Big \Vert u_{k}\Big \Vert 
_{KR}$ for frames/bases, we treat the question in Section 6 below for 
the classical spaces $ L^{2}(I^{d})$.\ \par 

\medskip

%Now, we extend this result to measure spaces over (almost) 
%arbitrary 1-dimensional "smooth manifold" of finite length, as follows.\ \par 
%\ \par 
\bf Proposition 5.1. \it Let $ \varphi :I\longrightarrow X$ 
be a continuous injection of $ I=[0,1]$ in a normed space $ X$ differentiable 
a.e. (with respect to Lebesgue measure $ dx$), and the distance on 
$ I$ \rm be \it defined by\ \par 
\ \par 
  \centerline{ $ \rho (x,y)=\Vert \varphi (x)-\varphi (y)\Vert _{X}$, 
$ x,y\in I$.}
\ \par 
\it Let further, $ \mu $ be a continuous (without point masses) probability 
measure on $ I$, satisfying\ \par 
\ \par 
  \centerline{ $ \displaystyle \int _{I}d\mu (y)\displaystyle \int _{y}^{1}\Vert 
\varphi '(x)\Vert _{X}dx=:$ $ C^{2}(\mu ,\varphi )<$ $ \infty $\rm .}
\ \par 
\it Then, every Bessel sequence $ u=(u_{k})$ in $ L^{2}(\mu )=:$ 
$ L_{0}^{2}((I,\mu )$ \rm fulfills\ \par 
\ \par 
  \centerline{ $ \displaystyle \sum _{k}\Vert u_{k}\Vert ^{2}_{KR}\leq 
B^{2}C(\mu ,\varphi )<$ $ \infty $\bf , }
\ \par 
\it where $ B(u)>0$ comes from the Bessel condition\rm .\ \par 
\ \par 
\bf Proof. \rm Following the proof of Theorem 3.1\it (1) 
\rm and using that for $ f\in C^{\infty }$,\ \par 
\ \par 
  \centerline{ $ \Lip(f)\leq 1\Leftrightarrow \vert 
f(x)-f(y)\vert \leq \Vert \varphi (x)-\varphi (y)\Vert \Leftrightarrow 
\vert f'(x)\vert \leq \Vert \varphi '(x)\Vert _{X}$ 
($ x\in I$),}
\ \par 
\rm we obtain, for every $ h\in L_{0}^{2}(\mu )$ and $ J_{\mu }(h)(x):=
\int _{0}^{x}hd\mu $,\ \par 
\ \par 
  \centerline{ $\Vert h\Vert _{KR}=\sup\Big \{\displaystyle \int _{I}f\,hd\mu 
:f\in C^{\infty },\Lip(f)\leq 1\Big \}=$ }
\ \par 
  \centerline{ $ =\sup\Big \{\displaystyle \int _{I}f'J_{\mu }(h)dx:
\Big \vert f'(x)\Big \vert \leq\Vert\varphi '(x)\Vert _{X}\Big \}=
\displaystyle \int _{I}\Big \vert J_{\mu }(h)\Big \vert \cdot\Vert\varphi 
'(x)\Vert _{X}dx\leq $}
\ \par 
  \centerline{ $ \leq\Vert J_{\mu }(h)\Vert _{L^{2}  (I,vdx)}$\bf ,}
 \ \par 
where $ v(x)=\Vert \varphi '(x)\Vert _{X}$. A mapping $ 
Th:=J_{\mu }(h)$, $ Th(x):=\int _{I}k(x,y)h(y)d\mu $ 
acting as $ T:L^{2}(\mu )\longrightarrow L^{2}(I,vdx)$ is in the Hilbert-Schmidt 
class $ {\cS}_{2}$ if and only if\ \par 
\ \par 
  \centerline{ $\Vert T\Vert _{2}^{2}=\displaystyle \int \displaystyle \int 
_{I\times I}\Big \vert k(x,y)\Big \vert ^{2}d\mu (y)v(x)dx=\displaystyle \int 
_{0}^{1}d\mu (y)\displaystyle \int _{y}^{1}v(x)dx=:C^{2}(\mu ,\varphi 
)<\infty .$}
\ \par 
\rm If $ u=(u_{k})$ is Bessel (with $ \sum _{k}\vert (h,u_{k})\vert 
^{2}\leq B(u)^{2}\Vert h\Vert ^{2}$, $ \forall h\in L^{2}_{w}$), 
and the last condition is fulfilled, then $ u_{k}=Ae_{k}$ where 
$ (e_{k})$ is orthonormal and $ \Vert A\Vert \leq B(u)$, and 
hence
\ \par 
 $$
 \sum _{k}\Vert u_{k}\Vert ^{2}_{KR}\leq 
\sum _{k}\Vert (TA)e_{k}\Vert _{2}^{2}\leq 
$$
$$
\Vert TA\Vert _{2}^{2}\leq\Vert T\Vert _{2}^{2}\Vert 
A\Vert ^{2}\leq B^{2}(u)C^{2}(\mu ,\varphi ).
$$
\ \par 
\bf Remark. \rm In particular, the following (known?) 
formula appeared in the proof:\ \par 
\ \par 
  \centerline{ $\Vert h\Vert _{KR}=\displaystyle \int _{I}\Big \vert 
J_{\mu }(h)\Big \vert \cdot\Vert\varphi '(x)\Vert _{X}dx$;}
\ \par 
\rm see also comments below.\ \par 
\ \par 
\bf II. Examples of interpolation spaces appearing conspicuously in Theorem 
3.2. \rm Lemma 3 above suggests that all decreasing rates of $ \Vert u_{k}\Vert 
_{KR}$ can really occur, and so all cases of convergence/divergence 
of $ \sum _{k}\varphi (\Vert u_{k}\Vert _{KR})$ are different and non 
empty. The following partial cases are of interest.\ \par 
\ \par 
\bf (1) \rm The most known interpolation space between $ l^{2}$ and 
$ l^{\infty }$ is $ l^{p}$, $ 2<p<\infty $, which is included in Theorem 
3.2 with\ \par 
\ \par 
  \centerline{ $ r(t)=t^{1-{\frac{2}{p}} }$;}
\ \par 
\rm it serves for the case of power-like decreasing of $ b_{n}(\Lip_{1})$, 
or $ s_{n}(T)$ (if $ \Lip_{1}=T(B(L^{2}))$), and consequently 
of $ \Vert u_{n}\Vert _{KR}$:\ \par 
\ \par 
  \centerline{ $ \log{\frac{1}{s_{n}}} \approx \log(n)$, 
$ n\longrightarrow \infty $.}
\ \par 
\rm In particular, point (2) of Theorem 3.1 (where $ \Omega =$ 
$ I^{d}$, $ d\geq 2$) can be seen now as a partial case of Theorem 
3.2 since, in the hypotheses of 3.1(2), $ \Lip_{1}=TB(L^{\infty })\supset 
TB(L^{2})$ and $ T\in \bigcap _{p>d}{\cS}_{p}(L^{2}\longrightarrow L^{2})$ 
(which was already observed in the proof of Theorem 3.1).\ \par 
\ \par 
\bf (2) \rm The following spaces $ l^{\varphi }$ of slowly decreasing 
sequences $ (s_{n})$ are conjectured to appear as $ s$-numbers (or 
Bernstein $ n$-widths) of $ \Lip_{1}$ for partial cases of the triples 
$ \Omega ={\Bbb T}^{\infty }$, $ \rho =\rho _{\epsilon }$, 
$ m_{\infty }$ described in the proof of Lemma 3 above:\ \par 
\ \par 
- $ \displaystyle \sum _{n}s_{n}^{C\,\log\log{\frac{\displaystyle 1}{\displaystyle 
s_{n}  }} }<\infty $ (corresponding to $ \log{\frac{1}{s_{n}}} 
$ $ \approx $ $ {\frac{\log(n)}{\log\log(n)}} $; the case is included 
in Theorem 3.2 with\ \par 
\ \par 
  \centerline{ $ r(t)=t\cdot \exp\Big \{-{\frac{\displaystyle 1}{\displaystyle 
C}} \cdot {\frac{\displaystyle \log(t^{2})}{\displaystyle \log\log(t^{2})}} 
(1+o(1)\Big \}$, as $ t\longrightarrow \infty $}
\ \par 
\rm (follows from the known $ b^{-1}(y)={\frac{y}{\log(y)}} (1+o(1))$ 
for $ b(x)=x\cdot \log(x)$), which is eventually concave (since 
$ t\longmapsto r(t)=o(t)$ for $ t\longrightarrow \infty $ and 
lies in the Hardy fields, see \cite{Bou1976}, L'Appendice du Ch.V);\ \par 
\ \par 
- $ \displaystyle \sum _{n}s_{n}^{C(\log{\frac{\displaystyle 1}{\displaystyle 
s_{n}  }} )^{\alpha }}<\infty $, $ \alpha >1$ (corresponding 
to $ \log{\frac{1}{s_{n}}} $ $ \approx $ $ (\log(n))^{1/\alpha }$; the 
case is included in Theorem 3.2 with\ \par 
\ \par 
  \centerline{ $ r(t)=t\cdot exp\Big \{-({\frac{\displaystyle 1}{\displaystyle 
C}} \cdot \log(t^{2}))^{1/\alpha }\Big \}$, }
\ \par 
\rm which is eventually concave as $ t\longrightarrow \infty $ (by 
the same argument as above);\ \par 
\ \par 
- $ \displaystyle \sum _{n}e^{-{\frac{\displaystyle C}{\displaystyle 
s_{n}^{\beta }  }} }<\infty $, $ \beta >0$ (corresponding to 
$ \log{\frac{1}{s_{n}}} $ $ \approx $ $ (c+{\frac{1}{\beta }} \log\log(n))$; 
the case is included in Theorem 3.2 with\ \par 
\ \par 
  \centerline{ $ r(t)=Ct/(\log(t^{2}))^{1/\beta }$, }
\ \par 
\rm which is eventually concave as $ t\longrightarrow \infty $ (by 
the same argument as above).\ \par 
\ \par 
\bf III. In terms of the Bernstein $ n$-widths. \rm It is quite easy 
to see that a part of Theorem 3.2, namely point (2), is still true 
with a (slightly?) relaxed hypothesis: we replace the assumption that 
$ \Lip_{1}$ is of the form $ \Lip_{1}\supset T(B(L^{2}))$ for 
a compact $ T$ with a hypothesis that the optimal subspaces for Bernstein 
widths $ b_{n}(\Lip_{1})$ are ordered by inclusion (see Section 2 above 
for the definitions): $ H_{n}(\Lip_{1})\subset $ $ H_{n+1}(\Lip_{1})$, 
$ n=$ $ 1,2,...$ Namely, the following property holds.\ \par 
\ \par 
\bf Proposition 5.2. \it Let $ \Omega ,\rho ,m$ be 
a compact probability triple for which there exist Bernstein optimal 
subspaces $ H_{n}(\Lip_{1})\subset L^{2}(\Omega ,m)$ such that\ \par 
\ \par 
  \centerline{ $ H_{n}(\Lip_{1})\subset $ $ H_{n+1}(\Lip_{1})$\rm , $ 
n=$ $ 1,2,...$}
\ \par 
\it Then there exists an orthonormal sequence $ (u_{k})_{k\geq 0}\subset $ 
$ \Lip(\Omega )\subset L^{2}_{{\Bbb R}}(\Omega ,m)$, such that\ \par 
\ \par 
  \centerline{ $ \Vert u_{n}\Vert _{KR}\geq b_{n}(\Lip_{1})$, 
$ n=1,2,...$}
\ \par 
\bf Proof. \rm Let $ e_{1}\in H_{1}$, $ \Vert e_{1}\Vert 
_{2}=b_{1}$, and assume that $ e_{k}$, $ k\leq n$ are chosen 
in a way that $ e_{k}\in H_{n}$, $ e_{k}\perp e_{j}$ ($ k\not= j$) 
and $ \Vert e_{k}\Vert _{2}=b_{k}$. Since $ b_{n+1}B(H_{n+1})\subset 
\Lip_{1}$, there exists a vector $ e_{n+1}\in H_{n+1}\ominus H_{n}\subset 
\Lip(\Omega )$ with $ \Vert e_{n+1}\Vert _{2}=b_{n+1}$ 
(and hence, $ e_{n+1}\in \Lip_{1}$). For the constructed sequence $ (e_{n})$, 
we set\ \par 
\ \par 
  \centerline{ $ u_{n}=e_{n}/b_{n}$}
\ \par 
\rm and obtain an orthonormal sequence $ (u_{n})\subset \Lip(\Omega 
)$ such that $ \Lip(u_{n})\leq 1/b_{n}$, and hence $ \Vert u_{n}\Vert 
_{KR}\geq \int _{\Omega }u_{n}e_{n}dm=b_{n}(\Lip_{1})$. 
\ \par 
\ \par 
\bf IV. Remark: an ``uncertainty inequality" for $ \Vert u\Vert _{KR}$. 
\rm As it is already used several times (in particular in the proof 
of 5.2 above), for a smooth function $ u\in \Lip(\Omega )$ the following 
inequality holds\ \par 
\ \par 
  \centerline{ $ \Vert u\Vert _{KR}\Lip(u)\geq \Vert u\Vert _{2}^{2}$.}
\ \par 
\rm Indeed, $ \Vert u\Vert _{KR}\geq \int _{\Omega }u(u/\Lip(u))dm$. 
\ \par 
As a consequence, one can observe that for every normalized 
Bessel sequence $ (u_{k})$, its $ \Lip$ norms must be sufficiently large, 
so that $ \sum _{k}\varphi ({\frac{1}{\Lip(u_{k})}} )<\infty $ 
for any monotone increasing function $ \varphi \geq 0$ for which $ \sum 
_{k}\varphi (\Vert u_{k}\Vert )<$ $ \infty $ (compare with the statements 
of Section 3).\ \par 
\ \par 
\bf V. Remark: an explicit formula for $ \Vert u\Vert _{KR}$. \rm There 
are some cases where the norm $ \Vert \cdot \Vert _{KR}$ can be explicitly 
expressed in term of the triple $ \Omega ,\rho ,m$. In particular, 
if $ \Lip_{1}=T(B(L^{\infty }(\Omega ,m))$ then\ \par 
\ \par 
  \centerline{ $\Vert u\Vert _{KR}=\Vert T^{*}u\Vert 
_{L^{1}  (\Omega ,m)}$, $ \forall u\in L^{1}(\Omega ,m)$.}
\ \par 
\rm Indeed, 
$$
\Vert u\Vert _{KR}=\sup\Big \{ \int _{\Om}ufdm:f\in 
\Lip_{1}\Big \}=\Vert T^{*}u\Vert _{L^{1}  (\Omega ,m)}.
$$

 In particular, such a formula holds for $ (\Omega ,m)=(I^{d},m_{d})$, 
as it is mentioned in the proof of Theorem 3.1 (the corresponding $ 
T(\sum _{k\not= 0}c_{k}e^{i(k,x)})=\sum _{k\not= 0}\vert k\vert 
c_{k}e^{i(k,x)}$ is a multiplier on $ L^{p}_{0}$); for $ d=1$, the 
formula is mentioned in \cite{Ver2004}.\ \par 
\ \par 
\bf VI. Yet another characteristic of a compact set. \rm The following 
compactness measure seems to be closely related to the estimates of 
$ \Vert u_{n}\Vert _{KR}$:\ \par 
\ \par 
  \centerline{ $ t(n)=\sup\Big \{r>0:\exists x_{j}\in \Lip_{1},
x_{i}\perp x_{k}(i\not= k),\Vert x_{j}\Vert \geq 
r,1\leq j\leq n\Big \}$, $ n\geq 1$\sc .}
\ \par 
\rm It is easy to see that $ \sqrt{n} b_{n}(\Lip_{1})\geq t(n)\geq 
b_{n}(\Lip_{1})$, and in principle, we can use $ t(n)$ instead 
of $ b_{n}$ in the proof of Proposition 5.2. We can also derive the 
existence of finite orthonormal sequences $ (e_{j})_{j=1}^{n}\subset 
\Lip(\Omega )$ such that $ \sum _{j=1}^{n}\varphi (\Vert e_{j}\Vert _{KR})\geq 
n\varphi (b_{n}(\Lip_{1}))$, $ n=1,2,...$\ \par

%%%%%%%%%%%%%%%%%%%%%%%%%%%%%%%%
%%%%%%%%%%%%%%%%%%%%%%%%%%%%%%%%

 \section{A summary, and the best $ KR$-norms behavior 
for frames/bases in $ L^{2}(I^{d})$.}
\label{thebest}

\bf (A) A summary of the worst (generic) behavior of the $ KR$-norms 
\rm (all these claims are already proved above). For every Bessel sequence 
$ (u_{k})$ in $ L^{2}(I^{d})$, we have for $ d=1$: $ \displaystyle \sum _{k}\Big \Vert 
u_{k}\Big \Vert _{KR}^{2}<$ $ \infty $, and for $ d>1$: $ \displaystyle \sum _{k}\Big \Vert 
u_{k}\Big \Vert _{KR}^{d+\epsilon }<$ $ \infty $, $ \forall \epsilon >0$.\ \par 
\, \, \it These claims are sharp: \rm for every compact triple 
$ {\rm (}\Omega {\rm ,}\rho {\rm ,m)}$ and for every sequence $ {\rm (}\epsilon 
_{{\rm k}}{\rm )}_{{\rm k}\geq {\rm 1}}$, $ \epsilon _{{\rm k}}\geq {\rm 0}$, 
such that $ \sum _{{\rm k}}\epsilon _{{\rm k}}^{{\rm 2}}{\rm <\, }\infty 
$, there exists an orthonormal sequence $ (u_{{\rm k}}{\rm )}_{{\rm k}\geq 
{\rm 1}}$ in $ {\rm L}^{{\rm 2}}_{{\Bbb R}}{\rm (}\Omega {\rm ,m)}$ 
such that $ \Big \Vert u_{k}\Big \Vert _{KR}\geq $ $ c\epsilon _{k},$ 
$ k=1,2,...$ ($ c>0$), and in $ L^{2}(I^{d})$ there exists an orthonormal 
sequence $ (u_{k})$ such that $ \displaystyle \sum _{k}\Big \Vert u_{k}\Big \Vert 
_{KR}^{d}=$ $ \infty $.\ \par 
\, \, For a generic compact triple $ \Omega {\rm ,}\rho {\rm ,m}$, 
we can only claim $ \lim_{k}\Big \Vert u_{k}\Big \Vert _{KR}=$ $ 0$ 
for every Bessel sequence in $ {\rm L}^{{\rm 2}}_{{\Bbb R}}{\rm (}\Omega 
{\rm ,m)}$. The property is sharp in the following sense: for every 
sequence $ {\rm (}\epsilon _{{\rm k}}{\rm )}_{{\rm k}\geq {\rm 1}}$, 
$ \epsilon _{{\rm k}}{\rm >0}$, with $ \lim _{{\rm k}}\epsilon _{{\rm k}}{\rm =\, 
0}$, \it there exists a compact triple $ {\rm (}\Omega {\rm ,}\rho {\rm ,m)}$ 
\rm (with usual properties) and an orthonormal sequence $ (u_{{\rm k}}{ )}_{{\rm k}\geq 
{\rm 1}}$ in $ {\rm L}^{{\rm 2}}_{{\Bbb R}}{\rm (}\Omega {\rm ,m)}$ 
such that $ \Big \Vert u_{k}\Big \Vert _{KR}=$ $ c\epsilon _{k},$ $ 
k=1,2,...$ ($ {\frac{1}{2\sqrt{2} }} \leq c\leq {\frac{\displaystyle 2\sqrt{\displaystyle 
2} }{\displaystyle \pi }} $).\ \par 
\ \par 

\bf (B) Bases/frames with the least possible $ KR$-norms. \rm For the 
best possible behavior of $ \Big \Vert u_{k}\Big \Vert _{KR}$ we replace the words ``for every Bessel sequence'' by the words ``there exists Bessel sequence'', meaning that we look for the fastest rate of decrease of $\{\Big \Vert  u_{k}\Big \Vert _{KR}\}$. Then  for 
bases/frames/Bessel sequences on $ L^{2}(I^{d})$, we have different 
summation properties, and for $d=1$ the threshold is $2/3$ (and not $2$ as above), as follows.\ \par 
\ \par 
\, \, {\bf Theorem 6.1. }{\it Let $ d=1,2,...$ and $ \alpha =$ 
$ {\frac{2d}{d+2}} $ ($ \alpha <2$). Then, (1) there exists an orthonormal 
basis $ (u_{k})$ in $ L^{2}(I^{d})$ such that $ \displaystyle \sum _{k}\Big \Vert 
u_{k}\Big \Vert _{KR}^{\alpha +\epsilon }<$ $ \infty $, $ \forall \epsilon 
>0$, but (2) $ \displaystyle \sum _{k}\Big \Vert u_{k}\Big \Vert _{KR}^{\alpha 
}=$ $ \infty $, for every frame $ (u_{k})$ in $ L^{2}(I^{d})$ (in particular, 
for every Riesz basis).}
\ \par\ \par

 Let $ (u_{n})$ be the Haar basis 
in $ L_{0}^{2}(I^{d})$ enumerated with the following notation:

\medskip

  \centerline{ $ h=\, \chi _{(0,1/2)\, }-\, \chi _{(1/2,1)}$}
  
\medskip
\noindent stands for the \it Haar basic wavelet \rm on $ I\subset {\Bbb R}$; 
taking a subset $ \sigma \subset \, D:=\, \{1,2,...,d\}$, 
$ \sigma \not= \, \emptyset$, and a multiindex $ k=\, (k_{1},k_{2},...,k_{d})\in 
{\Bbb Z}_{+}^{d}$, where $ 0\leq k_{s}<\, 2^{j}$ for every $ s$ 
and $ j\in {\Bbb Z}_{+}$, define \it the Haar functions $ (u_{n}):=$ 
$ (h_{j,k,\sigma })$ \rm as\ \par 
\ \par 
  \centerline{ $ h_{j,k,\sigma }(x)=\, 2^{dj/2}\displaystyle \prod _{s\in 
\sigma }h(2^{j}x_{s}-k_{s})\displaystyle \prod _{s\in D\backslash \sigma 
}\chi _{(0,1)}(2^{j}x_{s}-k_{s})$,}
\ \par 
\rm where $ x=$ $ (x_{1},x_{2},..,x_{d})\in I^{d}$. Then (see for example, 
\cite{Me1992}, Section 3.9), $ (u_{n})$ forms an orthonormal basis in $ L_{0}^{2}(I^{d})$ 
($ j$ and $ k$ run over all mentioned above values, $ \sigma $ runs 
a finite set of $ 2^{d}-1$ elements). Obviously,\ \par 
\ \par 
  \centerline{ $ \supp(h_{j,k,\sigma })=\, Q_{j,k}:=\, \{x\in 
{\Bbb R}^{d}:\, 2^{j}x-k\in I^{d}\}=\, \prod _{s=1}^{d}[k_{s}2^{-j},(k_{s}+1)2^{-j}]$.}
\ \par 
\, \bf Lemma. \rm Let $ u\in L^{\infty }(I^{d})$, $ \supp(u)\subset Q_{j,k}$ 
and $ \int _{I^{d}}udx=\, 0$. Then,\ \par 
\ \par 
  \centerline{ $ \Big \Vert u\Big \Vert _{KR}\leq \, {\frac{d}{2}} 
\Vert u\Vert _{\infty }2^{-(d+1)j}$.}
\ \par 

{\bf Proof}
\rm Since $ \int _{I^{d}}udx=$ $ 0$, 
we can restrict ourselves in the formula\ \par 
\ \par 
  \centerline{ $ \Big \Vert u\Big \Vert _{KR}=$ $ \sup\Big \{\displaystyle \int _{I}ufdx:$ 
$ Lip(f)\leq $ $ 1\Big \}$}
\ \par 
\rm to the functions $ f$ with $ f(l)=\, 0$, $ Lip(f)\leq 1$ where 
$ l=\, (k_{s}2^{-j})_{s=1}^{d}$, and so $ \vert f(x)\vert \leq \, 
\vert l-x\vert $, $ x\in Q_{j,k}$. Changing variables, we have\ \par 
\ \par 
  \centerline{ $ \Big \Vert u\Big \Vert _{KR}\leq $ $ \displaystyle \int _{Q_{j,0}}\Big \Vert 
u\Big \Vert _{\infty }\Big \vert x\Big \vert dx\leq \, \displaystyle \int _{Q_{j,0}}\Big \Vert 
u\Big \Vert _{\infty }\displaystyle \sum _{s=1}^{d}x_{s}dx=\, $ 
}
\ \par 
  \centerline{ $ =\, \Big \Vert u\Big \Vert _{\infty }{\frac{\displaystyle d}{\displaystyle 
2}} 2^{-2j}2^{-j(d-1)}=\, \Big \Vert u\Big \Vert _{\infty }{\frac{\displaystyle d}{\displaystyle 
2}} 2^{-j(d+1)}$\rm . }
%\end{proof} 

{\bf  Proof of Theorem 6.1}
\ \par 
\, \, \rm(1) Applying Lemma to $ u=\, h_{j,k,\sigma }$,\ \par 
\ \par 
  \centerline{ $ \Big \Vert h_{j,k,\sigma }\Big \Vert _{KR}\leq \, 2^{jd/2}{\frac{\displaystyle 
d}{\displaystyle 2}} 2^{-j(d+1)}$.}
\ \par 
\rm Summing up (with a $ \gamma >\alpha $, $ \alpha =$ $ {\frac{2d}{d+2}} 
$), we get\ \par 
\ \par 
  \centerline{ $ \displaystyle \sum _{n}\Big \Vert u_{n}\Big \Vert _{KR}^{\gamma 
}\leq \, \displaystyle \sum _{\sigma }\displaystyle \sum _{j\geq 0}\displaystyle \sum 
_{k}\Big \Vert h_{j,k,\sigma }\Big \Vert _{KR}^{\gamma }\leq \, \displaystyle \sum 
_{\sigma }\displaystyle \sum _{j\geq 0}2^{jd}\Big (2^{jd/2}{\frac{\displaystyle d}{\displaystyle 
2}} 2^{-j(d+1)}\Big )^{\gamma }<\, \infty $\bf .}

\medskip

\rm (2) Recall that the space $ L^{1}_{0}(I^{d})$ endowed 
with the $ KR$-norm is isometrically embedded into the dual space $ 
(Lip_{0})^{*}$ (with respect to the standard duality $ (u,f)=\, \int 
_{I^{d}}ufdm$).\ \par 
\, \, \, The plan of the proof (suggested by E. Gluskin) 
is the following: consider some metric properties of the embedding\ \par 
\ \par 
  \centerline{ $ E^{*}:L^{2}_{0}(I^{d})\longrightarrow \, (Lip_{0})^{*}$}
\ \par 
\rm and its predual embedding\ \par 
\ \par 
  \centerline{ $ E:Lip_{0\, }\longrightarrow \, L^{2}_{0}(I^{d})$}
\ \par 
\rm from two different points of view. Namely, assuming that there exists 
a frame $ (u_{k})$ in $ L_{0}^{2}(I^{d})$ such that $ \displaystyle \sum _{k}\Big \Vert 
u_{k}\Big \Vert _{KR}^{\alpha }<\, \infty $, we show that\ \par 
\ \par 
\bf (I) \rm embeddings $ E$, $ E^{*}$ are $ 2$-nuclear operators (see 
below) and the $ 2$-nuclear approximation numbers $ a^{(2)}_{N}(E^{*})$ 
decrease as $ o(1/N^{1/d})$ when $ N\longrightarrow \infty $;\ \par 
\bf (II) \rm on the other hand, one can see that - at least for $ N=\, 
2^{jd},\, j=\, 1,2,...$ - the numbers $ a^{(2)}_{N}(E)$ (which 
coincide with $ a^{(2)}_{N}(E^{*})$) cannot be less than $ cN^{-1/d}$.\ \par 
\, \, The above contradiction shows property (2) of Theorem 
6.1.\ \par 
\ \par 
\, \, \bf Proof of point (I). 
\rm A linear operator $ T:X\longrightarrow 
Y$ between Banach spaces $ X$ and $ Y$ is said $ p$\it -nuclear \rm if 
$ Tx=\, \sum _{k}T_{k}x$, $ x\in X$ (weak convergence), $ rank(T_{k})\leq 
\, 1$ and $ \sum _{k}\Vert T_{k}\Vert ^{p}<\, \infty $; $ 
\inf\Big \{\Big (\displaystyle \sum _{k}\Big \Vert T_{k}\Big \Vert ^{p}\Big )^{1/p}:\, 
over\, all\, such\, representations\Big \}=:\, \Big \Vert 
T\Big \Vert _{N(p)}$ is called its $ p$-norm. $ N$\it -th $ p$-nuclear 
approximation number of $ T$ \rm ($ N=\, 1,2,...$) is\ \par 
\ \par 
  \centerline{ $ a^{(p)}_{N}(T):=\, \inf\Big \{\Big \Vert T-A\Big \Vert _{N(p)}:\, 
A:X\longrightarrow Y,\, rank(A)<\, N\Big \}$.}
\ \par 
\, \, \rm Assume now that there exists a frame $ (u_{k})$ 
in $ L_{0}^{2}(I^{d})$ such that $ \displaystyle \sum _{k}\Big \Vert u_{k}\Big \Vert 
_{KR}^{\alpha }<$ $ \infty $ where $ \alpha =$ $ {\frac{2d}{d+2}} $. 
Let $ Sf=\, \sum _{k}(f,u_{k})u_{k}$ be the frame operator on 
$ L_{0}^{2}(I^{d})$; $ S$ is an isomorphism $ S:L_{0}^{2}(I^{d})\longrightarrow L_{0}^{2}(I^{d})$, 
and $ E^{*}S:\, L_{0}^{2}(I^{d})\longrightarrow \, (Lip_{0})^{*}$ 
is a $ 2$-nuclear operator,\ \par 
\ \par 
  \centerline{ $ E^{*}Sf=$ $ \sum _{k\geq 1}(f,u_{k})E^{*}u_{k}$,}
\ \par 
\rm since $ \Vert E^{*}u_{k}\Vert_{(Lip_0)^*} =\, \Vert u_{k}\Vert _{KR}$ 
and $ \alpha <2$. Moreover, letting $ (u_{k})$ in the decreasing order 
of $ \Big \Vert u_{k}\Big \Vert _{KR}$, we get $ \Big \Vert u_{k}\Big \Vert _{KR}^{\alpha 
}=\, o(1/k)$ (as $ k\longrightarrow \infty $), and hence\ \par 
\ \par 
  \centerline{ $ a^{(2)}_{N}(E^{*}S)^{2}\leq \, \displaystyle \sum _{k\geq 
N}\Big \Vert u_{k}\Big \Vert _{KR}^{2}\leq \, \Big \Vert u_{N}\Big \Vert 
_{KR}^{2-\alpha }\displaystyle \sum _{k\geq N}\Big \Vert u_{k}\Big \Vert _{KR}^{\alpha 
}=\, o({\frac{\displaystyle 1}{\displaystyle N^{2/\alpha -1}}} )$,}
\ \par 
\rm and $ a^{(2)}_{N}(E^{*}S)=$ $ o({\frac{\displaystyle 1}{\displaystyle N^{1/\alpha 
-1/2}}} )=\, o({\frac{\displaystyle 1}{\displaystyle N^{1/d}}} )$, 
as $ N\longrightarrow \infty $ and $1/\alpha= 1/2+ 1/d$. Since $ S$ is invertible, and $ \Big \Vert UTV\Big \Vert 
_{N(p)}\leq \, \Vert U\Vert \cdot \Vert T\Vert _{N(p)}\cdot \Vert V\Vert 
$ for every $ T,U,V$, we have\ \par 
\ \par 
  \centerline{ $ a^{(2)}_{N}(E^{*})=$ $ o({\frac{\displaystyle 1}{\displaystyle N^{1/d}}} 
)$, as $ N\longrightarrow \infty $.}

\bigskip

\, \, \bf Proof of point (II). 
\rm (The proof was suggested 
by E. Gluskin). We need to show that there exists a constant $ c>0$ 
such that for every operator $ A_{N}:\, Lip_{0}\longrightarrow \, 
L_{0}^{2}(I^{d})$, $ rank(A_{N})<\, N=\, 2^{jd}$ ($ j=1,2,...$), 
one has $ \Vert E-A_{N}\Vert _{N(2)}\geq \, cN^{-1/d}$. To this 
end, we construct two linear mappings $ V=V_{N}:{\Bbb R}^{N}\longrightarrow Lip_{0}$ 
and $ U=U_{N}:L_{0}^{2}(I^{d})\longrightarrow {\Bbb R}^{N}$ such that\ \par 
\ \par 
  \centerline{ $ UEV=\, id_{{\Bbb R}^{N}}$, $ \Vert V:{\Bbb R}^{N}\longrightarrow 
Lip_{0}\Vert \leq \, CN^{{\frac{1}{2}} +{\frac{1}{d}} }$, $ \Vert U:L_{0}^{2}(I^{d})\longrightarrow 
{\Bbb R}^{N}\Vert =\, 1$,}
\ \par 
\rm where $ C>0$ does not depend on $ N$.\ \par 
\, \, Having these mappings at hand, we get $ U_{2N}(E-A_{N})V_{2N}=$ 
$ id_{{\Bbb R}^{2N}}-\, B_{N}$, where $ rank(B_{N})<\, N$ 
and so\ \par 
\ \par 
  \centerline{ $ \Vert U_{2N}(E-A_{N})V_{2N}\Vert _{N(2)}=$ $ \Vert id_{{\Bbb R}^{2N}}-$ 
$ B_{N}\Vert _{N(2)}\geq \, N^{1/2}$,}
\ \par 
\rm and on the other hand,\ \par 
\ \par 
\, $ \Vert U_{2N}(E-A_{N})V_{2N}\Vert _{N(2)}\leq \Vert U_{2N}\Vert 
\cdot \Vert E-A_{N}\Vert _{N(2)}\Vert V_{2N}\Vert \leq \, C(2N)^{{\frac{1}{2}} 
+{\frac{1}{d}} }\Vert E-A_{N}\Vert _{N(2)}$,  which gives $ \Vert E-A_{N}\Vert _{N(2)}\geq \, cN^{-1/d}$.\ \par 
\ \par 
\, \, \bf Construction of the mappings $ V=V_{N}:{\Bbb R}^{N}\longrightarrow 
Lip_{0}$ and $ U=U_{N}:L_{0}^{2}(I^{d})\longrightarrow {\Bbb R}^{N}$, 
$ N=\, 2^{jd}$, $ j=1,2,...$. \rm We use the similar scaling procedure 
as in the above proof of part (1) of Theorem 6.1: let $ \psi $ be a 
smooth function on $ {\Bbb R}^{d}$ such that $ \supp(\psi )\subset \, 
Q_{0}=\, I^{d}$, $ \Vert \psi \Vert _{L^{2}  (I^{d}  )}=\, 1$, 
$ \int _{I^{d}}\psi dm=\, 0$, and, for every $ j\in {\Bbb Z}_{+}$,\ \par 
\ \par 
  \centerline{ $ \psi _{k}=\, \psi _{j,k}(x):=\, 2^{jd/2}\psi 
(2^{j}x-k)$, $ k\in K_{j}$, }
\ \par 
\noindent where $ K_{j}=\, \{k=(k_{1},...,k_{d})\in {\Bbb Z}_{+}^{d}$: 
$ 0\leq k_{s}<$ $ 2^{j}$ ($ 1\leq s\leq d$)$ \}$. Then, $ \psi _{k}$ 
($ k\in K_{j}$) have pairwise disjoint supports and form an orthonormal 
family in $ L_{0}^{2}(I^{d})$, $ card(K_{j})=\, 2^{jd}:=\, N$. 
Now, setting\ \par 
\ \par 
  \centerline{ $ Va=\, \sum _{k\in K_{j}}a_{k}\psi _{k}$, $ a\in 
{\Bbb R}^{N}$,}
\ \par 
\noindent we obtain\ \par 
\ \par 
  \centerline{ $ \Big \Vert Va\Big \Vert _{Lip}\leq c\cdot \sup_{x\in I^{d}}\Big \vert 
\nabla (Va)(x)\Big \vert =\, c\cdot \max_{k\in K_{j}}\sup_{x\in I^{d}}\Big \vert 
a_{k}\nabla \psi _{k}(x)\Big \vert \leq \, C2^{jd/2}2^{j}\Big \Vert a\Big \Vert 
_{{\Bbb R}^{N}}$,}
\ \par
 \noindent where $ c>0,\, C>0$ depend only on $ d$ (and the choice of 
$ \psi $), which gives the needed $ \Vert V:{\Bbb R}^{N}\longrightarrow Lip_{0}\Vert 
\leq $ $ CN^{{\frac{1}{2}} +{\frac{1}{d}} }$.\ \par 
\, \, For $ U=U_{N}:L_{0}^{2}(I^{d})\longrightarrow {\Bbb R}^{N}$, 
we let $ Uf=\, ((f,\psi _{k}))_{k\in K_{j}}$, and obviously get 
$ UEV=$ $ id_{{\Bbb R}^{N}}$ and $ \Vert U:L_{0}^{2}(I^{d})\longrightarrow {\Bbb R}^{N}\Vert 
=$ $ 1$.

%%%%%%%%%%%%%%%%%%%

\end{document}